\documentclass[12pt,reqno]{amsart}
\usepackage[colorlinks=true,
pdfstartview=FitV, linkcolor=cyan, citecolor=magenta,
urlcolor=blue]{hyperref}
\usepackage{amsmath,amsfonts,latexsym,amssymb}
\usepackage{mathrsfs}
\usepackage[latin1]{inputenc}
\usepackage[T1]{fontenc}
\usepackage{ae,aecompl}
\usepackage{braket}
\usepackage{comment}
\usepackage{color}
\usepackage{graphicx}
\usepackage{subfig}	
\newtheorem{theorem}{Theorem}[section]
\newtheorem{lemma}[theorem]{Lemma}
\newtheorem{proposition}[theorem]{Proposition}
\newtheorem{corollary}[theorem]{Corollary}

\newtheorem*{thm*}{\protect\theoremname}
\theoremstyle{plain}
\newtheorem*{cor*}{\protect\corollaryname}

\renewcommand{\geq}{\geqslant}

\long\def\@savemarbox#1#2{\global\setbox#1\vtop{\hsize\marginparwidth 
  \@parboxrestore\tiny\raggedright #2}}
\newcommand{\tr}{\operatorname{Tr}}

\newcommand{\TT}{\mathsf{T}}

\newcommand{\Real}{\mathbb R}

\newcommand{\PSL}{\mathsf{PSL}(d,\Real)}
\newcommand{\SL}{\mathsf{SL}(d,\Real)}

\newcommand{\ms}{\mathsf}

\newcommand{\g}{\gamma}

\providecommand{\corollaryname}{Corollary}
\providecommand{\theoremname}{Theorem}

\begin{document}

\title{Pressure metrics for cusped Hitchin components}
\author[Bray]{Harrison Bray}
\address{George Mason University, Fairfax, Virginia 22030, hbray@gmu.edu}
\author[Canary]{Richard Canary}
\address{University of Michigan, Ann Arbor, MI 41809, canary@umich.edu}
\author[Kao]{Lien-Yung Kao}
\address{George Washington University, Washington, D.C. 20052, lkao@gwu.edu}
\author[Martone]{Giuseppe Martone}
\address{University of Michigan, Ann Arbor, MI 41809, martone@umich.edu}
\thanks{Canary was partially supported by grant  DMS-1906441 from the National Science Foundation
and grant 674990 from the Simons Foundation. Parts of this work were done during the third author's visit to National Center for Theoretical Sciences (NCTS) in Taiwan. 
He would like to thank NCTS for their support and hospitality.}

\begin{abstract}
We study the cusped Hitchin component consisting of (conjugacy classes of) cusped Hitchin representations of a torsion-free geometrically finite
Fuchsian group $\Gamma$ into $\mathsf{PSL}(d,\mathbb R)$. We produce Riemannian pressure metrics associated to the first fundamental weight
and the first simple root. We produce a pressure path metric associated to the Hilbert length and describe its degeneracy. 
\end{abstract}


\maketitle

\setcounter{tocdepth}{1}
\tableofcontents
\section{Introduction}

In this paper, we construct pressure metrics on the cusped Hitchin component of Hitchin representations of a torsion-free Fuchsian lattice into $\mathsf{PSL}(d,\mathbb R)$.
The first two metrics are mapping class group invariant, analytic Riemannian metrics. These metrics are associated to the first fundamental
weight and the first simple root. Our third pressure metric is 
based on the Hilbert length. It is a mapping class group invariant path metric which is an analytic
Riemannian metric off of the self-dual locus. These constructions are based on earlier constructions of Bridgeman, Canary, Labourie and Sambarino
\cite{BCLS,BCLS2,BCS} in the case of Hitchin components of closed surface groups.

The main new technical difficulties involve the fact that 
while the geodesic flow of a closed hyperbolic surface may be coded by a finite Markov shift, there is no finite Markov coding of 
the geodesic flow of a geometrically finite hyperbolic surface.
Stadlbauer \cite{stadlbauer} and Ledrappier-Sarig \cite{ledrappier-sarig} 
provide a countable Markov coding of the (recurrent portion of the) geodesic flow of a finite area hyperbolic surface.
In a previous paper, we used these codings, work of Canary-Zhang-Zimmer \cite{CZZ} on cusped Hitchin representations, and the Thermodynamic Formalism
for countable Markov shifts,  to establish counting and equidistribution results  for cusped Hitchin representations. 
In this paper, we apply the theory developed in that paper to construct our pressure metrics.

The long-term goal of this project is to realize these
metrics as the induced metric on the strata at infinity of the metric completion of the Hitchin component of a closed surface group with its
pressure metric. In the classical setting, when $d=2$, Masur \cite{masur-wp} showed that the metric completion of Teichm\"uller space  of a closed surface $S$, with
the Weil-Petersson metric,  is the augmented Teichm\"uller space. The strata at infinity in the augmented Teichm\"uller space come from Teichm\"uller space of,
possibly disconnected, surfaces  obtained from pinching $S$ along a multicurve.  We hope that the Hilbert length pressure metric when $d=3$ may be more natural to study given its connection
to Hilbert geometry. When $d=3$, the Hitchin component of a closed surface is the space of
holonomy maps of convex projective structures on the surface. The strata at infinity of the augmented Hitchin component would then be cusped Hitchin components consisting of finite area convex projective structures obtained from pinching the surface along a multicurve. We hope to eventually establish an analogue of Masur's result in the higher rank setting.
(See \cite{canary-hitchin} for a more detailed description of the  conjectural geometric picture of the augmented Hitchin component.)
\smallskip

We now discuss our results more precisely.
We recall that if $\Gamma$ is a torsion-free, geometrically finite Fuchsian group (i.e. a discrete non-abelian finitely generated subgroup of $\mathsf{PSL}(2,\mathbb R)$), then a Hitchin representation is a representation
$\rho:\Gamma\to \mathsf{PSL}(d,\mathbb R)$ which admits a positive equivariant limit map $\xi:\Lambda(\Gamma)\to \mathcal F_d$ where $\Lambda(\Gamma)\subset\partial\mathbb H^2$ is the limit set of $\Gamma$ and
$\mathcal F_d$ is the space of $d$-dimensional flags. As in the closed case, they all arise as type-preserving deformations of the restriction of
an irreducible representation of $\mathsf{PSL}(2,\mathbb R)$ into $\mathsf{PSL}(d,\mathbb R)$.

The {\em Hitchin component} $\mathcal H_d(\Gamma)$ is the space of conjugacy classes of Hitchin representations of
$\Gamma$ into  $\mathsf{PSL}(d,\mathbb R)$. Fock and Goncharov, see the discussion in \cite[Sec 1.8]{fock-goncharov}, show that the Hitchin component is topologically a cell.
(When $d=3$, $\mathcal H_3(\Gamma)$ is parameterized by Marquis \cite{marquis-module}, when $\Gamma$ is a lattice, and more generally by Loftin and Zhang \cite{loftin-zhang}.
Bonahon-Dreyer \cite[Thm. 2]{bonahon-dreyer} and Zhang \cite[Prop. 3.5]{zhang-internal} explicitly describe variations of the Fock-Goncharov parametrization when
$\Gamma$ is cocompact, and their analyses should extend to our setting.) More generally, if $\mathsf{G}$ is a real-split Lie subgroup of $\mathsf{PSL}(d,\mathbb R)$, let
$\mathcal H(\Gamma, \mathsf{G})$ be the space of Hitchin representations with image in $\mathsf{G}$. (In particular, $\mathcal H_d(\Gamma)=\mathcal H(\Gamma,\mathsf{PSL}(d,\mathbb R))$ in this notation.) Fock-Goncharov \cite{fock-goncharov} and Hitchin \cite{hitchin} (see also \cite[\S 9.3]{GRW_positive}) show that $\mathcal H(\Gamma, \mathsf{G})$ is topologically a cell.

\begin{theorem}\label{FGcell}
If $\Gamma\subset\mathsf{PSL}(2,\mathbb R)$ is torsion-free and geometrically finite and $\mathsf{G}$ is a real-split Lie subgroup of $\mathsf{PSL}(d,\mathbb R)$, then 
the cusped Hitchin component $\mathcal H(\Gamma,\mathsf{G})$ is an analytic manifold diffeomorphic to $\mathbb R^m$ for some $m\in\mathbb N$.
\end{theorem}

If $$\mathfrak a=\left\{\vec x\in\mathbb R^d\mid \sum x_i=0\right\}$$
is the standard Cartan algebra for $\mathsf{PSL}(d,\mathbb R)$, let 
$$\Delta=\left\{ \phi=\sum_{i=1}^{d-1} a_i\alpha_i\ |\ a_i\ge 0\ \forall\ i,\ \sum a_i>0\right\}\subset \mathfrak{a}^*$$
where $\alpha_i$ is the simple root given by $\alpha_i(\vec x)=x_i-x_{i-1}$. 
Notice that $\Delta$ is exactly the collection of linear functionals which are strictly positive on the interior of the Weyl chamber
$$\mathfrak a^+=\left\{\vec x\in\mathfrak a\mid  x_1\geq \dots \geq x_d\right\}.$$
Consider the {\em Jordan projection}
$\nu:\mathsf{PSL}(d,\mathbb R)\to \mathfrak a^+$ given by 
$$\nu(A)=(\log \lambda_1(A),\ldots,\log\lambda_d(A))$$
where $\lambda_1(A)\ge \cdots \ge\lambda_d(A)$ are the (ordered) moduli of generalized eigenvalues of $A$.

If $\phi\in\Delta$ and $\rho\in \mathcal H_d(\Gamma)$, denote by $\ell^\phi_\rho(\gamma)=\phi(\nu(\rho(\gamma)))$ the $\phi$-length of $\gamma\in \Gamma$. We may define
the $\phi$-entropy of $\rho$ as
$$h^\phi(\rho)=\lim_{T\to\infty}\frac{\# R_T^\phi(\rho)}{T}$$
where $[\Gamma_{hyp}]$ is the set of conjugacy classes of hyperbolic elements in $\Gamma$, and
$$R_T^\phi(\rho)=\left\{[\gamma]\in[\Gamma_{hyp}]\ |\ \ell^\phi_\rho(\gamma)\le T\right\}.$$
Moreover, if
$\rho,\eta\in\mathcal H_d(\Gamma)$, we may define the $\phi$-pressure intersection
$$I^{\phi}(\rho,\eta)=\lim_{T\to\infty} \frac{1}{|R_T^{\phi}(\rho)|}\sum_{[\gamma]\in R_T^{\phi}(\rho)}
\frac{\ell_\eta^{\phi}(\gamma)}{\ell_{\rho}^{\phi}(\gamma)},$$
and a renormalized $\phi$-pressure intersection
$$J^{\phi}(\rho,\eta)=\frac{h^{\phi}(\eta)}{h^{\phi}(\rho)}I^{\phi}(\rho,\eta).$$

Our key tool in the construction of the pressure metric will be results of Bray, Canary, Kao and Martone \cite{BCKM} and Canary, Zhang and Zimmer \cite{CZZ} which
combine to prove that all these quantities vary analytically. See \cite[Section 8.1]{BCLS} for the analogous statement when $\Gamma$ is cocompact.

\begin{theorem}
\label{everything analytic}
If $\Gamma\subset\mathsf{PSL}(2,\mathbb R)$ is torsion-free and geometrically finite and
$\phi\in\Delta$, then $h^\phi(\rho)$ varies analytically over $\mathcal H_d(\Gamma)$ and
$I^\phi$ and $J^\phi$ vary analytically over $\mathcal H_d(\Gamma)\times\mathcal H_d(\Gamma)$.
Moreover, if $\rho,\eta\in\mathcal H_d(\Gamma)$, then
$$J^\phi(\rho,\eta)\ge 1$$ 
and $J^\phi(\rho,\eta)=1$ if and only if $\ell^\phi_\rho(\gamma)=\frac{h^{\phi}(\eta)}{h^{\phi}(\rho)}\ell^\phi_\eta(\gamma)$ for all $\gamma\in\Gamma$.
\end{theorem}

Given $\phi\in\Delta$, we define a pressure form on the Hitchin component, by letting
$$\mathbb P^{\phi}\big|_{T_{\rho}\mathcal H_d(\Gamma)}={\rm Hess}(J^{\phi}(\rho,\cdot)).$$
Since $J^\phi$ achieves its minimum along the diagonal, ${\mathbb P}^{\phi}$ will always be non-negative.
However, it will not always be non-degenerate. Typically, the most difficult portion of the proof of
the construction of a pressure metric is to verify non-degeneracy, or, more generally, to characterize which vectors
are degenerate.

We first consider the first fundamental weight $\omega_1\in\Delta$, given by $\omega_1(\vec x)=x_1$.
As a consequence of a much more general result, Bridgeman, Canary, Labourie and Sambarino \cite{BCLS} prove that $\mathbb P^{\omega_1}$ 
is non-degenerate on the Hitchin component of a convex cocompact Fuchsian group. 
We recall that the mapping class group $\mathrm{Mod}(\Gamma)$ is the group of
(isotopy classes of) orientation-preserving self-homeomorphisms of $\mathbb H^2/\Gamma$.

\begin{theorem}
\label{spectral nondegenerate}
If $\Gamma\subset\mathsf{PSL}(2,\mathbb R)$ is torsion-free and geometrically finite, then the pressure form ${\mathbb P}^{\omega_1}$ is non-degenerate, so
gives rise to a mapping class group invariant, analytic Riemannian metric on $\mathcal H_d(\Gamma)$.
\end{theorem}

Bridgeman, Canary, Labourie and Sambarino \cite{BCLS2} later expanded their techniques to show that the first simple root gives rise
to a non-degenerate pressure metric on the Hitchin component of a closed surface group. We implement their outline in the cusped setting.
We make crucial use of a result of Canary, Zhang and Zimmer \cite{CZZ2} which assures us that simple root entropies are constant on the
Hitchin components of Fuchsian lattices (which generalizes a result of Potrie and Sambarino \cite{potrie-sambarino} for Hitchin components of closed surface groups).

\begin{theorem}
\label{simple nondegenerate}
If $\Gamma\subset\mathsf{PSL}(2,\mathbb R)$ is a torsion-free lattice, then the pressure form ${\mathbb P}^{\alpha_1}$ is non-degenerate, so
gives rise to a mapping class group invariant, analytic Riemannian metric on $\mathcal H_d(\Gamma)$.
\end{theorem}

Finally, we consider the functional $\omega_H$ associated to the Hilbert length given by $\omega_H(\vec x)=x_1-x_d$. 
It is easy to see that  if $C:\mathcal H_d(\Gamma)\to\mathcal H_d(\Gamma)$
is the contragredient involution and $\vec v\in T\mathcal H_d(\Gamma)$ is {\em anti-self-dual}, i.e. $\mathrm DC(\vec v)=-\vec v$, 
then $\mathbb P^{\omega_H}(\vec v,\vec v)=0$ (see \cite[Lem. 5.22]{BCS}).
In particular, $\mathbb P^{\omega_H}$ is not globally non-degenerate. However, one can still show that the pressure form
gives rise to a path metric. Bridgeman, Canary and Sambarino \cite[Sec. 5.8]{BCS} previously remarked that this is the
case when $\Gamma$ is a closed surface group.

\begin{theorem}
\label{hilbert mostly nondegenerate}
If $\Gamma\subset\mathsf{PSL}(2,\mathbb R)$ is torsion-free and geometrically finite, then
$\mathbb P^{\omega_H}$ gives rise to a mapping class group invariant path metric on $\mathcal H_d(\Gamma)$ which is an analytic Riemannian metric off of the self-dual locus.
\end{theorem}

When $d=3$, cusped Hitchin representations of a torsion-free lattice are holonomy maps of finite area convex projective surfaces and the Hilbert length is the translation length with respect to
the Hilbert metric. In this case, the analogy with the augmented Teichm\"uller space is most compelling and we expect that this case may be the easiest case in which to begin
the analysis of the augmented Hitchin component. Notice that our proposed augmented Hitchin component would be a proper subspace of the augmented Hitchin component 
introduced and studied in \cite{loftin-zhang}.

Theorems \ref{spectral nondegenerate} and \ref{hilbert mostly nondegenerate} are derived by generalizing the main result
of \cite[Thm. 1.4]{BCLS} into the setting of (cusped) Anosov representations of geometrically finite Fuchsian groups.  A (cusped)  $P_{1,d-1}$-Anosov representation of a 
geometrically finite Fuchsian group $\Gamma$ is type-preserving, i.e. takes hyperbolic elements to biproximal elements and parabolic elements to (weakly) unipotent elements, and admits an equivariant limit map from the limit set of $\Gamma$ into the partial flag variety whose elements are pairs $(L,H)$  where $L$ is a  line contained in a hyperplane 
$H$ in $\mathbb R^d$. See Section \ref{AnosovGF} for detailed definitions.

Let $\widetilde{\mathcal P}^{irr}_{\{1,d-1\}}(\Gamma,d)$ be the set of irreducible  $P_{\{1,d-1\}}$-Anosov representations into $\mathsf{PSL}(d,\mathbb R)$ and let
${\mathcal P}^{irr}_{\{1,d-1\}}(\Gamma,d)=\widetilde{\mathcal P}^{irr}_{\{1,d-1\}}(\Gamma,d)/\mathsf{PSL}(d,\mathbb R)$. If $\mathsf{H}$ is a reductive subgroup of $\mathsf{PSL}(d,\mathbb R)$, then an
element $h\in\mathsf{H}$ is $\mathsf{H}$-generic if its centralizer is a maximal torus in $\mathsf{H}$. If $\mathsf{H}=\mathsf{PSL}(d,\mathbb R)$, then an element is $\mathsf{H}$-generic
if and only if it is diagonalizable over $\mathbb C$ with distinct eigenvalues. A representation into $\mathsf{H}$ is said to be $\mathsf{H}$-generic if its image contains an $\mathsf{H}$-generic element.
In particular, all Hitchin representations are $\mathsf{PSL}(d,\mathbb R)$-generic, so Theorem \ref{spectral nondegenerate}  is a special case of the following
more general result.

\begin{theorem}
\label{BCLS generalized}
Suppose that $\Gamma\subset\mathsf{PSL}(2,\mathbb R)$ is torsion-free and geometrically finite.
If $W$ is an analytic submanifold of $\mathcal P^{irr}_{\{1,d-1\}}(\Gamma,d)$, $\mathsf H$ is a reductive subgroup of $\mathsf{PSL}(d,\mathbb R)$  and  every representation in
$W$ has image in $\mathsf{H}$ and is $\mathsf{H}$-generic, then $\mathbb P^{\omega_1}|_{\TT W}$ is an analytic Riemannian metric on $W$. Moreover, if $W$ is invariant under
a subgroup $M$ of the mapping class group, then  $\mathbb P^{\omega_1}|_{\TT W}$ is $M$-invariant.
\end{theorem}

The proof of Theorem \ref{BCLS generalized} follows the same outline as the proof of the main result in \cite{BCLS}.
Standard results from the Thermodynamic Formalism imply that $\vec v\in T_\rho W$ has 
\hbox{$\omega_1$-pressure} norm zero if and only if there exists $K$ so
that $D_{\vec v}\ell^{\omega_1}_\gamma=K\ell^{\omega_1}_\gamma$ for all hyperbolic elements $\gamma\in\Gamma$. We can then apply \cite[Lemma 9.8]{BCLS} to
show that $K=0$. An analysis using Labourie's cross ratio functions is used to show that this implies that $\vec v$ itself must
be zero, which completes the proof.

Finally, we remark that if $\Gamma$ is geometrically finite but has torsion, then it has a finite index normal subgroup $\Gamma_0$ which is torsion-free.
One may identify $\Gamma/\Gamma_0$ with a finite index subgroup $G$ of the mapping class group of $\mathbb H^2/\Gamma_0$ and then identify
$\mathcal H_d(\Gamma)$ with the submanifold of $\mathcal H_d(\Gamma_0)$ which is stabilized by $G$. It follows that
one obtains mapping class group invariant analytic Riemannian metrics $\mathbb P^{\omega_1}$ and $\mathbb P^{\alpha_1}$ on
$\mathcal H_d(\Gamma)$ and a mapping class group invariant path metric on $\mathcal H_d(\Gamma)$ which is analytic Riemannian
off of the self-dual locus.

\subsection*{Historical remarks}
Thurston described a metric on Teichm\"uller space which was the ``Hessian of the length of a random geodesic.'' Wolpert \cite{wolpert} showed that this
metric gives a scalar multiple of the classical Weil-Petersson metric. Bonahon \cite{bonahon} reinterpreted Thurston's metric in terms of geodesic currents.
McMullen \cite{mcmullen-wp} showed that one may interpret Thurston's metric in terms of Thermodynamic Formalism, as the Hessian of a pressure
intersection function. Bridgeman \cite{bridgeman-wp} generalized McMullen's construction to the setting of quasifuchsian space. Bridgeman, Canary, Labourie
and Sambarino \cite{BCLS} then showed how to use his construction to produce analytic Riemannian metrics at ``generic'' smooth points of deformation\
spaces of projective Anosov representations, and in particular on Hitchin components. Pollicott and Sharp \cite{pollicott-sharp-wp} gave an alternate
interpretation of this metric.

Kao \cite{kao-pm} used countable Markov codings to construct pressure metrics on Teichm\"uller spaces of punctured surfaces. Bray, Canary and Kao \cite{BCK}
generalized this to the setting of cusped quasifuchsian groups.

\subsection*{Acknowledgement} We thank the referee for their helpful comments, which allowed us to simplify our proofs and improve our results.

\section{Background}

\subsection{Linear algebra}\label{linear algebra} The {\em Jordan projection} $\nu\colon \SL\to\frak a^+$ is the map which associates to $
A\in\SL$ the list $(\log \lambda_1(A),\dots, \log \lambda_d(A))$ of logarithms of moduli of generalized eigenvalues of $A$ in decreasing order.

The {\em Cartan projection} $\kappa\colon\SL\to\frak a^+$ is
\[
\kappa(A)=\left(\log\sigma_1(A), \ldots,\log \sigma_d(A)\right)
\] 
where $\{\sigma_i(A)\}_{i=1}^d$ are the singular values of $A$ labelled in decreasing order. 
Recall that each element of $\mathsf{SL}(d,\mathbb R)$ may be written as $A=KDL$ where $K,L\in\mathsf{SO}(d)$ and $D$ is the diagonal matrix with
$(i,i)$-entry given by $\sigma_i(A)$. If $\alpha_k(\kappa(A))>0$, then $U_k(A)=K(\langle e_1,\ldots, e_k\rangle)$ is well-defined and is the $k$-plane spanned by
the first $k$ major axes of $A(S^{d-1})$.

Suppose that $\theta$ is a symmetric subset of $\{1,\ldots,d-1\}$, i.e. $k\in\theta$ if and only if $d-k\in\theta$. Define the {\em $\theta$-Cartan subspace} as
$$\mathfrak{a}_\theta=\{\vec a\in\mathfrak{a} : \alpha_j(\vec a)=0\text{ if }j\notin\theta\}$$
and let $\Delta_\theta$ denote the set of functionals in $\mathfrak{a}_\theta^*$ which are positive on the interior $\mathfrak a_\theta^+=\mathfrak a^+\cap \mathfrak a_\theta$.
In particular,
$$\Delta=\Delta_{\{1,\ldots,d\} }.$$

The {\em $\theta$-Cartan projection} $\kappa_\theta\colon\SL\to\mathfrak{a}_\theta$ is the unique  map
so that $\omega_k(\kappa_\theta(A))=\omega_k(\kappa(A))$ for all $A\in\mathsf{SL}(d,\mathbb R)$ and all $k\in\theta$.

If $\theta=\{k_1,\ldots, k_n\}$ we define the {\em $\theta$-flag variety}
$$
\mathcal F_\theta=\{ (F^{k_1}, F^{k_2},\dots, F^{k_n}) : F^{k_1}\subset F^{k_2}\subset\cdots\subset F^{k_n}\}$$
where each $F^{k_i}$ is a vector subspace of $\mathbb R^d$ of dimension $k_i$. In particular, the full flag variety $\mathcal F_d$
is the same as $\mathcal F_{\{1,2,\ldots,d-1\}}$ in this notation. 

Quint \cite{quint-ps} introduced a vector valued smooth cocycle, called the {\em $\theta$-Iwasawa cocycle},
$$B_\theta:\mathsf{SL}(d,\mathbb R)\times\mathcal F_\theta\to \mathfrak{a}_\theta$$
with the defining property that if $k\in\theta$, $A\in\mathsf{SL}(d,\mathbb R)$, $F\in\mathcal F_\theta$,
$\vec v_k$ is a non-trivial vector in $E^k(F^{k})\subset E^k(\mathbb R^d)$, where $E^k$ denotes the $k^{\rm th}$ exterior power, then
$$\omega_k(B_\theta(A,F))=\log\frac{||E^kA(\vec v_k)||}{||\vec v_k||}$$
where $||\cdot||$ is the norm on $E^k\mathbb R^d$ induced by the standard Euclidean norm on $\mathbb R^d$.
Note that the Jordan and Cartan projections (resp. $\theta$-Iwasawa cocycle) descend to well-defined functions on $\PSL$ (resp. $\PSL\times\mathcal F_\theta$).

\subsection{Thermodynamic Formalism}
In this section, we recall the background results we will need from the Thermodynamic Formalism for countable Markov shifts
as developed by Gurevich-Savchenko \cite{gur-sav}, Mauldin-Urbanski \cite{MU} and Sarig \cite{sarig-2009}.

Given a countable alphabet $\mathcal A$ and a transition matrix $\mathbb T=(t_{ab})\in\{0,1\}^{\mathcal A\times\mathcal A}$
a one-sided Markov shift is
$$\Sigma^+=\{x=(x_i)\in\mathcal A^{\mathbb N}\ |\ t_{x_ix_{i+1}}=1\ {\rm for}\ {\rm all}\ i\in\mathbb N\}$$
equipped with a shift map $\sigma:\Sigma^+\to\Sigma^+$ which takes $(x_i)_{i\in\mathbb N}$ to
$(x_{i+1})_{i\in\mathbb N}$. 
One says that $(\Sigma^+,\sigma)$ is {\em topologically mixing} if for all $a,b\in \mathcal A$, there exists $N=N(a,b)$
so that if $n\ge N$, then there exists $x\in\Sigma^+$ so that $x_1=a$ and $x_n=b$. The shift $(\Sigma^+,\sigma)$ has the
big images and pre-images property (BIP) if there exists a finite subset $\mathcal B\subset\mathcal A$ so
that if $a\in\mathcal A$, then there exists $b_0,b_1\in\mathcal B$ so that $t_{b_0,a}=1=t_{a,b_1}$.

Given a one-sided countable Markov shift $(\Sigma^+,\sigma)$ and a function  $g:\Sigma^+\to \mathbb R$, 
we say that $g$ is {\em locally H\"older continuous} if there exists $C>0$ and $\eta\in (0,1)$
so that if $x,y\in\Sigma^+$ and $x_i=y_i$ for all $ 1\le i\le n$, then
$$|g(x)-g(y)|\le C\eta^n.$$ 
If $n\in\mathbb N$, the {\em $n^{th}$-ergodic sum}  of $g$ at $x\in\Sigma^+$ is
$$S_ng(x)=\sum_{i=1}^{n} g(\sigma^{i-1}(x))$$
and $\mathrm{Fix}^n=\{x\in\Sigma^+\ |\ \sigma^n(x)=x\}$ is the set of periodic words with period dividing $n$.

The {\em pressure} of a 
locally H\"older continuous function $g:\Sigma^+\to\mathbb R$ is defined to be
$$P(g)=\sup \left\{ h_\sigma(m)+\int_{\Sigma^+} g \ dm : m\in \mathcal{M}_\sigma\ \textrm{and}\  -\int_{\Sigma^+} g \ dm<\infty\right\}$$ where $\mathcal M_\sigma$ is the space of $\sigma$-invariant probability measures on $\Sigma^+$ and
$h_\sigma(m)$ is the measure-theoretic entropy of $\sigma$ with respect to the measure $m$.

A  $\sigma$-invariant Borel probability measure $m$ on $\Sigma^+$ is an {\em equilibrium measure} for 
a locally H\"older continuous function $g:\Sigma^+\to\mathbb R$ if 
$$P(g)=h_\sigma(m)+\int_{\Sigma^+} g \ dm.$$ We remark that there are several different but equivalent definitions of pressure and equilibrium measure in the current setting. Readers can find a more detailed discussion of this in Bray-Canary-Kao-Martone \cite[pg. 11]{BCKM}.
Mauldin-Urbanski (\cite[Thm. 2.6.12,\ Prop. 2.6.13\ and\  2.6.14]{MU}) and
Sarig (\cite[Cor. 4]{sarig-2003}, \cite[Thm 5.10\ and\  5.13]{sarig-2009}) prove that the pressure function is
real analytic in our setting and compute its derivatives.
Recall that $\{g_u:\Sigma^+\to \mathbb R\}_{u\in M}$ is a {\em real analytic family} if $M$ is a real analytic
manifold and for all $x\in \Sigma^+$, $u\to g_u(x)$ is a real analytic function on $M$. 

\begin{theorem}[Mauldin-Urbanski, Sarig]
\label{pressure analytic}
Suppose that  $(\Sigma^+,\sigma)$ is a one-sided countable Markov shift which has (BIP) and is topologically mixing.
If $\{g_u:\Sigma^+\to \mathbb R\}_{u\in M}$ is a real analytic family of locally H\"older continuous functions such
that $P(g_u)<\infty$ for all $u$, then $u\to P(g_u)$ is real analytic.

Moreover, if $\vec v\in T_{u_0}M$ and there exists a neighborhood $U$ of $u_0$ in $M$ so that if $u\in U$,
then $-\int_{\Sigma^+} g_u dm_{g_{u_0}}<\infty$, then
$$D_{\vec v}P(g_u)=\int_{\Sigma^+} D_{\vec v}(g_u(x))\  dm_{g_{u_0}}.$$
\end{theorem}
In the case of finite Markov shifts, the assumption that $P(g_u)<\infty$ is automatically satisfied and Theorem \ref{pressure analytic} is due to Ruelle \cite{ruelle-book} and Parry-Pollicott \cite{parry-pollicott}.

\medskip

Bowen and Series \cite{bowen-series} constructed a finite Markov coding for the action of a convex cocompact group $\Gamma$  on
its limit set $\Lambda(\Gamma)$. Dal'bo and Peign\'e \cite{dalbo-peigne}, when $\Gamma$ is geometrically finite but not a lattice, and 
Stadlbauer \cite{stadlbauer} and Ledrappier-Sarig \cite{ledrappier-sarig}, when $\Gamma$ is a lattice,
constructed a countable Markov coding for the action of $\Gamma$ on its conical limit set $\Lambda_c(\Gamma)$.
We summarize their crucial properties below (see \cite{BCK} for a more complete description in our language). If $a\in\mathcal A$, then $G(a)$  is the associated 
element of $\Gamma$ and
$\log r(a)$ is  ``coarsely'' the translation distance (of some fixed basepoint) of $G(a)$.

\begin{theorem}[Bowen-Series \cite{bowen-series}, Dal'bo-Peign\'e \cite{dalbo-peigne}, Ledrappier-Sarig \cite{ledrappier-sarig}, Stadlbauer \cite{stadlbauer})]
\label{coding properties}
Suppose that $\Gamma$ is a torsion-free geometrically finite Fuchsian group.
There exists a topologically mixing Markov shift $(\Sigma^+,\sigma)$ with countable alphabet $\mathcal A$ with (BIP)  and maps 
$$G:\mathcal A\to\Gamma,\ \ \omega:\Sigma^+\to\Lambda(\Gamma),\ \  \mathrm{and}\ \ r:\mathcal A\to \mathbb N$$
with the following properties.
\begin{enumerate}
\item
$\omega$ is locally H\"older continuous, finite-to-one and $\omega(\Sigma^+)=\Lambda_c(\Gamma)$, i.e.  the complement in $\Lambda(\Gamma)$ of the set of fixed points
of parabolic elements of $\Gamma$. 
Moreover, $\omega(x)=G(x_1)\omega(\sigma(x))$ for every
$x\in\Sigma^+$.
\item
If $x\in \mathrm{Fix}^n$, then $\omega(x)$ is the attracting fixed point of $G(x_1)\cdots G(x_n)$. Moreover, if $\gamma\in\Gamma$ is hyperbolic,
then there exists  $x\in\mathrm{Fix}^n$ (for some $n$) so that $\gamma$ is conjugate to $G(x_1)\cdots G(x_n)$ and $x$ is unique up to cyclic permutation.
\item
There exists $D\in\mathbb N$ so that $1\le \#(r^{-1}(n))\le D$ for all $n\in\mathbb N$.
\end{enumerate}
\end{theorem}

\subsection{Anosov representations of geometrically finite Fuchsian groups}
\label{AnosovGF}
We next recall the definition of a $P_\theta$-Anosov representation of a geometrically finite Fuchsian group
and the results of  Bray-Canary-Kao-Martone \cite{BCKM}  and  Canary-Zhang-Zimmer \cite{CZZ} which will play a crucial role in our work.

Let $\Gamma$ be a geometrically finite Fuchsian group and let $\theta$ be a symmetric subset of $\{1,\ldots,d-1\}$.
We say that a representation $\rho:\Gamma\to\mathsf{PSL}(d,\mathbb R)$ is {\em $P_\theta$-Anosov}, if there exists a
continuous $\rho$-equivariant map $\xi_\rho:\Lambda(\Gamma)\to \mathcal F_\theta$ so that
 \begin{enumerate}
 \item
 $\xi_\rho$ is {\em transverse}, i.e. if $x\ne y\in\Lambda(\Gamma)$ and $k\in\theta$, then
 $$\xi^k_\rho(x)\oplus\xi^{d-k}_\rho(y)=\mathbb R^d,$$
\item
 $\xi_\rho$ is {\em strongly dynamics preserving}, i.e. if
  $\{\gamma_n\}$ is a sequence in $\Gamma$ so
 that $\gamma_n(b_0)\to x\in\Lambda(\Gamma)$ and $\gamma_n^{-1}(b_0)\to y\in\Lambda(\Gamma)$ for some basepoint $b_0\in\mathbb H^2$, then
 if $F\in\mathcal F_\theta$ is transverse to $\xi_\rho(y)$, then $\rho(\gamma_n)(F)\to\xi_\rho(x)$.
 \end{enumerate}
We denote the space of $P_\theta$-Anosov representations of $\Gamma$ into $\mathsf{PSL}(d,\mathbb R)$ by $\widetilde{\mathcal P}_\theta(\Gamma,d)$.
We will need the following observation, which follows immediately from the above definition.

\begin{lemma}
\label{schottky subgroup}
If $\rho:\Gamma\to \mathsf{PSL}(d,\mathbb  R)$ is in $\widetilde{\mathcal P}_{\{1,d-1\}}(\Gamma,d)$ and $\Gamma_0$ is a Schottky subgroup of $\Gamma$,
then $\rho|_{\Gamma_0}$ is a projective Anosov representation of the convex cocompact subgroup $\Gamma_0$.
\end{lemma}
Canary, Zhang and Zimmer establish fundamental properties of $P_\theta$-Anosov representations of  geometrically finite Fuchsian groups which generalize the properties of classical Anosov representations.

\begin{theorem}[Canary-Zhang-Zimmer \cite{CZZ}]
\label{cusped Hitchin properties}
Suppose that $\Gamma$  is a geometrically finite Fuchsian group,
$\rho:\Gamma\to\mathsf{PSL}(d,\mathbb R)$ is a $P_\theta$-Anosov representation.
\begin{enumerate}
\item
If $\gamma\in\Gamma$ is hyperbolic and $k\in\theta$, then $\rho(\gamma)$ is $P_k$-proximal.
\item
If $\alpha\in\Gamma$ is parabolic, then $\rho(\alpha)$ is weakly unipotent in $\mathsf{PSL}(d,\mathbb R)$, i.e.  its Jordan-Chevalley decomposition has elliptic semi-simple part and non-trivial unipotent part.
\item
There exist $A,a>0$ so that if $\gamma\in\Gamma$ and $k\in\theta$, then
$$Ae^{ad(b_0,\gamma(b_0))}\ge e^{\alpha_k(\kappa_\theta(\rho(\gamma)))}\ge \frac{1}{A}e^{\frac{d(b_0,\gamma(b_0))}{a}}$$
where $b_0$ is a basepoint for $\mathbb H^2$.
\item
$\rho$ has the \emph{$P_\theta$-Cartan property},  i.e.  whenever $\{\gamma_n\}$ is a sequence of distinct elements of $\Gamma$
such that $\gamma_n(b_0)$ converges to $z\in\Lambda(\Gamma)$, then $\xi^k_\rho(z)=\lim U_k(\rho(\gamma_n))$ for all $k\in\theta$.
\end{enumerate}
\end{theorem}

They also show that limit maps of Anosov representations vary analytically.

\begin{theorem}[Canary-Zhang-Zimmer \cite{CZZ}]
\label{limit map analytic}
If  $\{\rho_u:\Gamma\to\mathsf{PSL}(d,\mathbb R)\}_{u\in M}$ is a real analytic family of $P_\theta$-Anosov representations of a geometrically finite Fuchsian group and 
$z\in\Lambda(\Gamma)$, then the map from $M$ to  $\mathcal F_\theta$
given by $u\to \xi_{\rho_u}(z)$ is  real analytic.
\end{theorem}

If $\rho\in\widetilde{\mathcal P}_\theta(\Gamma,d)$, the $\theta$-{\em Benoist limit cone} of $\rho$ is
\[
\mathcal B_\theta(\rho)=\bigcap_{n\geq 0}\overline{\bigcup_{\|\kappa_\theta(\rho(\gamma))\|\geq n}\mathbb R_+\kappa_\theta(\rho(\gamma))}\subset \mathfrak a^+_\theta.
\]
The positive dual to the $\theta$-Benoist limit cone is given by
\begin{equation}\label{dual Benoist}
\mathcal B_\theta(\rho)^+=\left\{\phi\in\mathfrak a_\theta^*\mid \phi\left(\mathcal B_\theta(\rho)-\{\vec 0\}\right)\subset(0,\infty)\right\}.
\end{equation}

In previous work \cite{BCKM}, we constructed potentials on the Markov shift which encode the
spectral properties of Anosov representations of geometrically finite, torsion-free Fuchsian groups. 
First we use the $\theta$-Iwasawa cocycle to define a vector-valued roof function $\tau_\rho:\Sigma^+\to\mathfrak{a}_\theta$ by
$$\tau_\rho(x)=B_\theta\big(\rho(G(x_1)),\rho(G(x_1))^{-1}(\xi_\rho(\omega(x)))\big)$$
If $\phi\in\mathcal B_\theta(\rho)^+$ one defines the roof function $\tau_\rho^\phi=\phi\circ\tau_\rho$.
Notice that since $\mathcal B_\theta(\rho)$ is contained in the interior of the positive Weyl chamber $\mathfrak{a}_\theta^+$, the set
$\Delta_\theta$ is contained in $\mathcal B_\theta(\rho)^+$.

We use the Thermodynamic Formalism for countable Markov shifts to analyze these potentials. In particular, we use a renewal theorem
of Kesseb\"ohmer and Kombrink  \cite{kess-komb} to generalize arguments of Lalley \cite{lalley} to establish counting and equidistribution
results in our setting. We summarize the results we will need from our work below.

\begin{theorem}[Bray-Canary-Kao-Martone \cite{BCKM}]
\label{roof properties}
Suppose that $\Gamma$  is a torsion-free, geometrically finite Fuchsian group which is not convex cocompact,
$\rho:\Gamma\to\mathsf{PSL}(d,\mathbb R)$ is a $P_\theta$-Anosov representation 
and $\phi\in\mathcal B_\theta(\rho)^+$. Then, there exists a locally H\"older continuous function
$\tau_\rho^\phi=\phi\circ\tau_\rho:\Sigma^+\to\mathbb R$ such that
\begin{enumerate}
\item
$\tau_\rho^\phi$ is eventually positive, i.e. there exist $N\in\mathbb N$ and $B>0$ such that $S_n\tau^\phi_\rho(x)>B$ for all $n\geq N$ and $x\in\Sigma^+$.
\item
There exists $d(\phi)>0$, so that  the function $h\mapsto P(-h\tau_\rho^\phi)$ is finite, proper and strictly monotone on $(d(\phi),\infty)$ and infinite otherwise.
\item
There exists $C_\rho>0$, and for all $x_1\in\mathcal A$, $c(\rho,\phi,x_1)\ge 1/d(\phi)$ so that if $x\in\Sigma^+$, then
$$\Big|\tau_\rho^{\phi}(x)-c(\rho,\phi,x_1)\log r(x_1)\Big|\le C_\rho.$$
\item
If $x=\overline{x_1\cdots x_n}$ is a periodic element of $\Sigma^+$, then
$$S_n\tau_\rho^\phi(x)=\ell^\phi_\rho\big(G(x_1)\cdots G(x_n)\big).$$
\item
The $\phi$-entropy  $h^\phi(\rho)$ of $\rho$ is the unique solution of $P(-h\tau_\rho^\phi)=0$.
Moreover,  
$$\lim_{T\to\infty} \frac{h^\phi(\rho) T R_T^\phi(\rho)}{e^{h^\phi(\rho)T}}=1.$$
\item
There is a unique equilibrium measure $m_\rho^\phi$ for $-h^\phi(\rho)\tau_\rho^\phi$.
\end{enumerate}
\end{theorem}

We also established a rigidity theorem for renormalized pressure intersection and use our  equidistribution result to
give a thermodynamical reformulation of the pressure intersection. Recall that $\Delta_\theta\subset\mathfrak{a}_\theta^*$ is the set of functionals
which are strictly positive on the interior of $\mathfrak{a}_\theta^*$. We often restrict to this set of linear functionals, since
$\Delta_\theta$ is the largest set of linear functionals contained in $\mathcal B_\theta(\rho)^+$ for 
every $P_\theta$-Anosov representation $\rho$.

\begin{theorem}[Bray-Canary-Kao-Martone \cite{BCKM}]
\label{intersection rigidity}
If $\rho,\eta:\Gamma\to \mathsf{PSL}(d,\mathbb R)$ are $P_\theta$-Anosov representations of a geometrically finite
Fuchsian group and $\phi\in\Delta_\theta$, then
$$J^{\phi}(\rho,\eta)\ge 1$$
with equality if and only if 
$$\ell^{\phi}_\rho(\gamma)=\frac{h_\phi(\eta)}{h_\phi(\rho)}\ell^{\phi}_\eta(\gamma)$$
for all $\gamma\in\Gamma$. 
Moreover,
$$I^\phi(\rho,\eta)=\frac{\int_{\Sigma^+} \tau_\eta^\phi\ dm_\rho^\phi}{\int_{\Sigma^+} \tau_\rho^\phi\ dm_\rho^\phi}$$
and $-I^\phi(\rho,\eta)$ is the slope of the tangent line at $(h^\phi(\rho),0)$ to
\[
\mathcal C^\phi(\rho,\eta)=\{(a,b)\in\mathbb R^2\mid P(-a\tau_\rho^\phi-b\tau^\phi_\eta)=0, a\geq 0, b\geq 0, a+b>0\}.
\]
\end{theorem}

\subsection{Hitchin representations}
We say that a basis $b=(b_1,\ldots,b_d)$ is {\em consistent} with a pair $(F,G)$ of transverse flags if
$\langle b_i\rangle=F^i\cap G^{d-i+1}$ for all $i$. We denote by $U(b)_{>0}\subset \SL$ the subsemigroup of upper triangular unipotent matrices which are totally positive with respect to $b$, i.e. $A\in U(b)_{>0}$ if, in the basis $b$, $A$ is upper triangular unipotent and the determinants of all the minors of $A$ are positive, unless they are forced to be zero by the fact that $A$ is upper
triangular.

Then, a $k$-tuple $(F_1,\ldots,F_k)$ in $\mathcal F_d$ is  {\em positive} if there exists a basis $b$ consistent
with $(F_1,F_k)$ and there exists $\{u_2,\ldots,u_{k-1}\}\in U(b)_{>0}$ so that 
$F_i=u_{k-1}\cdots u_iF_k$ for all $i=2,\ldots, k-1$. If $X$ is a subset of $S^1$, we say that  a map $\xi:X\to\mathcal F_d$ is {\em positive} if whenever
$(x_1,\ldots,x_k)$ is a consistently ordered $k$-tuple in $X$ (ordered either clockwise or counter-clockwise),
then $(\xi(x_1),\ldots,\xi(x_k))$ is a positive $k$-tuple of flags.

Let $\Gamma$ be a geometrically finite Fuchsian group and let  $\Lambda(\Gamma)\subset\partial\mathbb H^2$ be its limit set.
Following Fock and Goncharov \cite{fock-goncharov}, a {\em Hitchin representation}  $\rho:\Gamma\to \mathsf{PSL}(d,\mathbb R)$
is a representation 
such that there exists a $\rho$-equivariant positive map $\xi_\rho:\Lambda(\Gamma)\to\mathcal F_d$.
If $S$  is closed, Hitchin representations are just the traditional Hitchin representations introduced by Hitchin \cite{hitchin} and further studied
by Labourie \cite{labourie-invent}. When $\Gamma$ contains a parabolic element, we sometimes refer to these Hitchin representations as cusped Hitchin
representations to distinguish them from the traditional Hitchin representations. 

Canary, Zhang and Zimmer \cite{CZZ} proved the following important structural results. (Sambarino~\cite{sambarino-positive} independently showed
that Hitchin representations are strongly irreducible.)

\begin{theorem}\label{thm:hitchinborel}
If $\rho\in\mathcal H_d(\Gamma)$, then $\rho$ is $\{1,\ldots,d-1\}$-Anosov and strongly irreducible.
\end{theorem}

We recall that Sambarino \cite[Theorem A]{sambarino-positive} classified the possible Zariski closures of images of Hitchin representations.

\begin{theorem}[Sambarino \cite{sambarino-positive}]\label{thm: Sambarino}
Suppose that $\Gamma\subset\mathsf{PSL}(2,\mathbb R)$ is a geometrically finite Fuchsian group, and $\rho:\Gamma\to\mathsf{PSL}(d,\mathbb R)$ is a Hitchin representation. 
Then the Zariski closure of $\rho(\Gamma)$ either  lies in an irreducible image of $\mathsf{PSL}(2,\mathbb R)$ or is
conjugate to either $\mathsf{PSL}(d,\mathbb R)$, $\mathsf{PSp}(2n,\mathbb R)$ when $d=2n$, $\mathsf{PSO}(n,n-1)$ when $d=2n-1$, or $\mathsf{G}_2$ when $d=7$.
\end{theorem}

\medskip\noindent
{\bf Historical remarks:}
The results in this subsection generalize earlier results in the case when $\Gamma$ is  convex cocompact.
More precisely, when $\Gamma$ is convex cocompact,  Theorem \ref{cusped Hitchin properties} follows from work of Labourie \cite{labourie-invent}, 
Fock and Goncharov \cite{fock-goncharov}, Guichard and Wienhard \cite{guichard-wienhard}, Kapovich, Leeb and Porti \cite{klp-anosov,klp-morse actions}, Gu\'erituad, Guichard, Kassel, and Wienhard \cite{ggkw} and 
Tsouvalas \cite{tsouvalas}, Theorems \ref{limit map analytic} and \ref{intersection rigidity} are due to Bridgeman, Canary, Labourie, and Sambarino \cite{BCLS} and 
Theorem \ref{roof properties} is due to Sambarino \cite{sambarino-quantitative}. Finally, Theorem \ref{thm:hitchinborel} is due to Labourie \cite{labourie-invent}.

Anosov representations of geometrically finite Fuchsian groups are also relatively Anosov in the sense of Kapovich and Leeb \cite{kapovich-leeb} and relatively dominated in the sense of Zhu \cite{zhu-reldom}.
In particular, one can derive Theorem \ref{cusped Hitchin properties} in either of their settings. The approach in \cite{CZZ} was motivated by the need
to prove Theorem \ref{limit map analytic},  whose proof was not clear from either pre-existing viewpoint. Zhu and Zimmer \cite{zhu-zimmerI} have now generalized the techniques
of \cite{CZZ} to establish a generalization of Theorem \ref{limit map analytic} to the setting of all relatively Anosov representations.

\section{Entropy, intersection and the pressure form}

Our pressure form is defined as the Hessian of a renormalized intersection function, so it is crucial to show that this function is analytic (or at least $C^2$).
Let $\widetilde{\mathcal P}_\theta(\Gamma,d)$ be the space of $P_\theta$-Anosov representations of $\Gamma$ into $\mathsf{PSL}(d,\mathbb R)$.

\begin{theorem}
If $\tilde W$ is an analytic submanifold of $\widetilde{\mathcal P}_\theta(\Gamma,d)$  and  $\phi\in\Delta_\theta$, then $h^\phi(\rho)$ varies analytically over $\tilde W$ and
$I^\phi$ and $J^\phi$ vary analytically over $\tilde W\times \tilde W$.
Moreover, if $\rho,\eta\in \tilde W$, then
$$J^\phi(\rho,\eta)\ge 1$$ 
and $J^\phi(\rho,\eta)=1$ if and only if $\ell^\phi_\rho(\gamma)=\frac{h^{\phi}(\eta)}{h^{\phi}(\rho)}\ell^\phi_\eta(\gamma)$ for all $\gamma\in\Gamma$.
\end{theorem}

\begin{proof} If $\Gamma$ is convex cocompact, then this result is established in \cite{BCLS}. So we will assume that $\Gamma$ is geometrically finite but not convex cocompact for the rest of this proof.

Theorem \ref{limit map analytic} implies that the
limit map $\xi_\rho$ varies analytically  over $\tilde W$.
Since $\tau_\rho(x)=B_\theta(\rho(G(x_1)),\rho(G(x_1))^{-1}(\xi_\rho(\omega(x))))$ and $B_\theta$ is analytic we see that $\tau_\rho$, and  hence 
$\tau_\rho^\phi=\phi\circ\tau_\rho$, varies analytically over
$\tilde W$. It then follows from Theorem \ref{roof properties} and Theorem \ref{pressure analytic} that $(h,\rho)\to P(-h\tau_\rho^\phi)$ is analytic on
$(d(\phi),\infty)\times\tilde W$. Since $P(-h^\phi(\rho)\tau^\phi_\rho)=0$ and
$$\frac{d}{dt} P(-t\tau_\rho^\phi)\Big|_{t=h^\phi(\rho)}=-\int \tau_\rho^\phi\ dm_\rho^\phi<0$$
for all $\rho\in\tilde W$, the Implicit Function Theorem implies that $h^\phi(\rho)$ varies analytically over $\tilde W$.

Let
$$R=\tilde W\times \tilde W\times\hat D_\phi\qquad\text{where}\qquad \hat D_\phi=\left\{(a,b)\in\mathbb R^2\ | a+b>d(\phi)\right\}$$
and let $P_R:R\to \mathbb R$ be given by $P_R(\rho,\eta,a,b)=P(-a\tau^\phi_\rho-b\tau^\phi_\eta)$. Mauldin and Urbanski \cite[Thm. 2.1.9]{MU} show that if $f$ is locally
H\"older continuous, then $P(f)$ is finite if and only if $Z_1(f)<+\infty$, where
$$Z_1(f)=\sum_{s\in\mathcal A}e^{\sup\{ f(x)\ :\ x_1=s\}}<+\infty.$$
By grouping the terms so that $r(s)=n$, Theorem \ref{roof properties} implies that 
$$Z_1(-a\tau^\phi_\rho-b\tau^\phi_\eta)\le \sum_{n=1}^\infty  e^{-(a+b)\left(\frac{1}{d(\phi)}\log n -\max\{C_\rho,C_\eta\}\right)}$$ 
so  $P(-a\tau^\phi_\rho-b\tau^\phi_\eta)<+\infty$ if $a+b>d(\phi)$. 
Therefore, $P_R$ is finite on $R$, and hence, by Theorem \ref{pressure analytic}, analytic on $R$.
As above, $P_R$ is a submersion on $P_R^{-1}(0)$, so $P_R^{-1}(0)$ is an analytic submanifold of $R$.
Moreover, by Theorem \ref{intersection rigidity},  $-I^\phi(\rho,\eta)$ is the slope of the tangent line to 
$P_R^{-1}(0)\cap \{(\rho,\eta)\times\hat D_\phi\}$ at the point $(\rho,\eta,(h(\rho),0))$, so $I^\phi(\rho,\eta)$ is analytic. Since entropy is analytic,
it follows immediately that $J^\phi(\rho,\eta)$ is analytic. 

The final claim follows directly from Theorem \ref{intersection rigidity}.
\end{proof}

Let $\widetilde{\mathcal P}^{irr}_\theta(\Gamma,d)$ be the set of irreducible representations in $\widetilde{\mathcal P}_\theta(\Gamma,d)$ and let
${\mathcal P}^{irr}_\theta(\Gamma,d)=\widetilde{\mathcal P}^{irr}_\theta(\Gamma,d)/\mathsf{PSL}(d,\mathbb R)$. The argument of \cite[Proposition 7.1]{BCLS} shows that the action of $\mathsf{PSL}(d,\mathbb R)$ on $\widetilde{\mathcal P}^{irr}_\theta(\Gamma,d)$ is free, proper and analytic. It follows that if $W$
is an analytic submanifold of ${\mathcal P}^{irr}_\theta(\Gamma,d)$, then its pre-image $\tilde W$ is an analytic submanifold of $\widetilde{\mathcal P}_\theta(\Gamma,d)$.

In this setting, we get the following result which generalizes Theorem \ref{everything analytic} from the introduction.

\begin{corollary}
\label{everything analytic theta}
If $W$ is an analytic submanifold of $\mathcal P^{irr}_\theta(\Gamma,d)$  and  $\phi\in\Delta_\theta$, then $h^\phi(\rho)$ varies analytically over $W$ and
$I^\phi$ and $J^\phi$ vary analytically over $W\times  W$.
Moreover, if $\rho,\eta\in W$, then
$$J^\phi(\rho,\eta)\ge 1$$ 
and $J^\phi(\rho,\eta)=1$ if and only if $\ell^\phi_\rho(\gamma)=\frac{h^{\phi}(\eta)}{h^{\phi}(\rho)}\ell^\phi_\eta(\gamma)$ for all $\gamma\in\Gamma$.
\end{corollary}

Given $\phi\in\Delta_\theta$, we define a pressure form on any analytic submanifold  $W$ of $\mathcal P^{irr}_\theta(\Gamma,d)$
by letting
$$\mathbb P^{\phi}\big|_{\TT_\rho  W}={\rm Hess}(J^{\phi}(\rho,\cdot)).$$
If $\vec v=\frac{d}{dt}\big\vert_{t=0} [\rho_t]$ where $\{\rho_t\}_{t\in (-\epsilon,\epsilon)}$ is a one-parameter analytic family in $W$,
then
\[
\mathbb P^\phi(\vec v,\vec v)=\frac{d^2}{dt^2}\Big|_{t=0} J^{\phi}(\rho_0,\rho_t).
\]

We note that the exact same definitions apply when $\Gamma$ is a cocompact lattice, see \cite[Sect. 5.5]{BCS}.
We observe the following immediate properties.

\begin{proposition}
\label{pressure form basics}
If $W$ is an analytic submanifold of $\mathcal P^{irr}_\theta(\Gamma,d)$  and  $\phi\in\Delta_\theta$,
then
$\mathbb P^\phi$ is analytic and non-negative, i.e. if $\vec v\in \TT W$,
then $\mathbb P^\phi(\vec v,\vec v)\ge 0$. Moreover, if $M$ is a subgroup of the mapping class group of $\Gamma$ and $W$ is $M$-invariant,
then $\mathbb P^\phi$ is $M$-invariant.
\end{proposition}

\begin{proof}
The pressure form
$\mathbb P^\phi$ is analytic, since $J^\phi$ is analytic and is non-negative since $J^\phi$ achieves its minimum along the diagonal, see Theorem \ref{intersection rigidity}.
If $\psi\in M$, then $\ell_{\gamma}^\phi(\rho\circ\psi)=\ell^\phi_{\psi^{-1}(\gamma)}(\rho)$,
so $R_T^\phi(\rho\circ\psi)=\psi^{-1}\big(R_T^\phi(\rho)\big)$ and $h^\phi(\rho)=h^\phi(\psi\circ\rho)$, so $J^\phi(\rho,\eta)=J^\phi(\rho\circ\psi,\eta\circ\psi)$, which implies
that $\mathbb P^\phi$ is $\psi$-invariant. 
\end{proof}

The following degeneracy criterion for pressure metrics is standard in the setting of finite Markov shifts, see for example \cite[Cor. 2.5]{BCS}, but
requires a little more effort in the setting of countable Markov shifts. In our setting, this criterion is established exactly as in  Lemma 8.1 in \cite{BCK}.

\begin{proposition}
\label{degeneracy condition}
Suppose that $W$ is an analytic submanifold of $\mathcal P^{irr}_\theta(\Gamma,d)$ and  $\phi\in\Delta_\theta$.
If $\vec v\in \TT W$ and $\phi\in\Delta_\theta$, then ${\mathbb P}^{\phi}(\vec v,\vec v)=0$ if and only if 
$${\rm D}_{\vec v} \left( h^{\phi}\ell^{\phi}_\gamma\right) =0$$
for all $\gamma\in\Gamma$.
\end{proposition}

We next observe that ${\rm D}_{\vec v}\log\ell^{\phi}_\gamma$ is independent of $\gamma$ if $\vec v$ is degenerate and $\gamma$ is hyperbolic.

\begin{lemma}
\label{improved degeneracy condition}
Suppose that $W$ is an analytic submanifold of $\mathcal P^{irr}_\theta(\Gamma,d)$  and  $\phi\in\Delta_\theta$.
If $\vec v\in \TT W$ and ${\mathbb P}^{\phi}(\vec v,\vec v)=0$, then, 
$${\rm D}_{\vec v} \ell^{\phi}_\gamma=K\ell_\gamma^{\phi}$$
for all $\gamma\in\Gamma$, where $K=-\frac{{\rm D}_{\vec v} h^\phi}{h^{\phi}(\rho)}$.
\end{lemma}

\begin{proof}
By Proposition \ref{degeneracy condition}, if $\vec v$ is degenerate, then
$$h^\phi(\rho){\rm D}_{\vec v}\ell_\gamma^{\phi}= -\big({\rm D}_{\vec v} h^\phi \big)\ell_\gamma^{\phi}$$
for all hyperbolic $\gamma\in\Gamma$, so
${\rm D}_{\vec v}\ell_\gamma^{\phi}= K\ell_\gamma^{\phi}$
for all hyperbolic  $\gamma\in\Gamma$. If $\gamma\in\Gamma$ is  parabolic, then $\ell_\gamma^\phi$ is the zero function, so the condition holds
trivially.
\end{proof}

Recall that if $M$ is a real analytic manifold, 
an analytic function $f:M \rightarrow \mathbb R$  has \hbox{{\em log-type} $K$} at $v\in\ms T_uM$ 
if $f(u)\ne 0$ and
$$
{\rm D}_u{\log(|f|)}(v)=K\log(|f(u)|).
$$
In this language, Lemma \ref{improved degeneracy condition} implies that 
if $\Lambda^{\phi}_\gamma$ is the function defined by $\rho\to e^{\ell_\rho^{\phi}(\gamma)}$ and ${\mathbb P}^\phi(\vec v,\vec v)=0$ then there exists $K$
so that $\Lambda^{\phi}_\gamma$ has log-type $K$ at $\vec v$ for all $\gamma\in\Gamma$.

\section{The spectral radius pressure metric}

We are now ready to establish our  generalization of the main theorem of \cite{BCLS}.

\medskip\noindent
{\bf Theorem \ref{BCLS generalized}.} {\em
If $W$ is an analytic submanifold of $\mathcal P_{1,d-1}^{irr}(\Gamma,d)$, $\mathsf H$  is a reductive subgroup of $\mathsf{PSL}(d,\mathbb R)$  and  every representation in
$W$ has image in $\mathsf{H}$ and is $\mathsf{H}$-generic, then $\mathbb P^{\omega_1}|_{\TT W}$ is an analytic Riemannian metric on $W$. Moreover, if $W$ is invariant under
a subgroup $M$ of the mapping class group, then  $\mathbb P^{\omega_1}|_{\TT W}$ is $M$-invariant.}

\medskip

Since every cusped Hitchin representation is irreducible and $\mathsf{SL}(d,\mathbb R)$-generic, Theorem \ref{spectral nondegenerate} follows immediately from Theorem \ref{BCLS generalized}.

\begin{proof}[Proof of Theorem \ref{BCLS generalized}]
Proposition \ref{pressure form basics} implies that
we only need to prove that every non-zero vector $\vec v\in \TT W$ is non-degenerate. 
Suppose $\vec v_0\in \TT_{\rho_0}W$ is a degenerate vector. Lemma \ref{improved degeneracy condition} implies that 
there exists $K$ so that $\Lambda^{\omega_1}_\gamma$ has log-type $K$ at $\vec v_0$ for all $\gamma\in\Gamma$. Our first step consists of showing that $K=0$.

Let $\rho_0(\alpha)$ be an $\mathsf{H}$-generic element. We claim that $\alpha$ is hyperbolic and, consequently that $\rho_0(\alpha)$ is biproximal. 
If not, $\alpha$ is parabolic, and $\rho_0(\alpha)$ is weakly unipotent, see Theorem \ref{cusped Hitchin properties} part (2). 
In particular,  $\rho_0(\alpha)$ is not diagonalizable over $\mathbb C$, so its centralizer cannot be a  maximal torus, which contradicts the assumption that $\rho_0(\alpha)$ is $\mathsf{H}$-generic.

Let $\beta$ be an element of $\Gamma$, so that $\alpha$ and $\beta$ generate a free convex cocompact subgroup $\Gamma_0$ of $\Gamma$.
Lemma \ref{schottky subgroup} then implies that the restriction $\rho_0|_{\Gamma_0}$ is $P_{\{1,d-1\}}$-Anosov in the traditional sense.
One may then choose an open neighborhood  $W_0$ of $\rho_0$ in $W$ so that if $\rho\in W_0$, then $\rho|_{\Gamma_0}$ is 
$P_{\{1,d-1\}}$-Anosov and $\rho(\alpha)$ is $\mathsf{H}$-generic.
We then can consider the analytic family $\{\rho|_{\Gamma_0}\}_{\rho\in W_0}$ and apply  \cite[Lemma 9.8]{BCLS} to  conclude that $K=0$.

Next, we show that if $\Lambda_\gamma^{\omega_1}$ is log-type zero at $\vec v_0$ for all $\gamma\in\Gamma$, then $\vec v_0=0$.

Recall that Labourie defines a continuous cross-ratio ${\bf b}$ on pairs of mutually transverse hyperplanes $(P_1,P_2)$ and lines $(L_1,L_2)$ by setting
$${\bf b}(P_1,P_2,L_1,L_2)=\frac{\phi_1(\vec v_1)\phi_2(\vec v_2)}{\phi_1(\vec v_2)\phi_2(\vec v_1)}$$
where $\phi_1$ and $\phi_2$ are linear functionals with kernel $P_1$ and $P_2$ and $\vec v_1$ and $\vec v_2$ are non-zero vectors in $L_1$ and $L_2$.
One may then
define a continuous cross-ratio ${\bf b}_\rho:\Lambda(\Gamma)^{(4)}\to \mathbb R$, where $\Lambda(\Gamma)^{(4)}$ is the set of pairwise distinct quadruples in $\Lambda(\Gamma)$, by setting
$${\bf b}_\rho(x,y,z,w)={\bf b}(\xi_\rho^{d-1}(x),\xi_\rho^{d-1}(y),\xi_\rho^1(z),\xi_\rho^1(w)).$$
Note that $\bf b_\rho$ is well-defined because $\rho$ is $P_{\{1,d-1\}}$-Anosov.
Suppose that $(\alpha,\beta)$ is a pair of hyperbolic elements of $\Gamma$ generating a rank two  Schottky subgroup of $\Gamma$,
then $(\rho(\alpha),\rho(\beta))$ generate a projective Anosov Schottky group, so \cite[Prop. 10.4]{BCLS} gives that 
$${\bf b}_\rho(\alpha^-,\beta^-,\beta^+,\alpha^+)
=\lim_{n\to\infty}\frac{\Lambda^{\omega_1}_{\alpha^n\beta}(\rho)}{\Lambda^{\omega_1}_{\alpha^n}(\rho)}.$$
It follows that ${\bf b}_\rho(\alpha^-,\beta^-,\beta^+,\alpha^+)$ is of log-type zero (since ratios of log-type zero functions are log-type zero, as are limits of log-type zero functions).
Since such quadruples are dense in $\Lambda(\Gamma)^{(4)}$, it follows that
${\bf b}_\rho(x,y,z,w)$ is log-type zero for all quadruples in $\Lambda(\Gamma)^{(4)}$.

Given a projective frame $F=(L_1,\ldots, L_{d+1})$ for $\mathbb P(\mathbb R^d)$ and a projective frame $F^*=(P_1,\ldots,P_{d+1})$ for $\mathbb P((\mathbb R^d)^*)$,
one can define a smooth injective immersion 
$$\mu_{F,F^*}:\mathsf{PSL}(d,\mathbb R)\to W(d)$$
where $W(d)$ is the quotient of the space of $(d+1)\times(d+1)$ matrices via a multiplicative action of $(\mathbb R-\{0\})^{2(d+1)}$ whose action on the coefficients of the matrix is 
given by $(a_1,\cdots,a_{d+1},b_1,\cdots,b_{d+1})(M_{ij})=a_ib_j M_{ij}$ (see \cite[Section 10.2]{BCLS}). 
(Recall that a projection frame for $\mathbb P(\mathbb R^d)$ is a collection of $d+1$ lines so that no $d$ lines are contained in any hyperplane.)
Specifically one chooses non-zero vectors $\vec v_i\in L_i$ and
covectors $\phi_i\in P_i$ so that $\sum \vec v_i=0$ and $\sum \phi_i=0$ and defines 
\[
\mu_{F,F^*}(A)=[\phi_i(A(\vec v_j))].
\] 
This smooth injective immersion and the cross ratio ${\bf b}_\rho$ are related by the following crucial property, whose proof proceeds exactly as in \cite[Lemma 10.7]{BCLS}. 

\begin{lemma}
Suppose $\{x_1,\ldots, x_{d+1}\}$ and $\{y_1,\ldots, y_{d+1}\}$ are two collections of pairwise distinct points in $\Lambda(\Gamma)$,  and that 
$F=\{\xi_\rho^1(x_1),\dots, \xi_\rho^{1}(x_{d+1})\}$ and $F^*=\{\xi_\rho^{d-1}(y_1),\dots, \xi_\rho^{d-1}(y_{d+1})\}$ are projective frames. Then
$$\mu_{F,F^*}(\rho(\alpha))=[{\bf b}_\rho(y_i,z,\alpha(x_j),w)]$$
for arbitrary $z,w\in\Lambda(\Gamma)$ and for all $\alpha\in\Gamma$.
\end{lemma}

The remainder of the proof then simply mimics the proof of \cite[Lemma 10.8]{BCLS}.

Let $\{\rho_t\}$ be a path in $W$ so that $\frac{d}{dt}\rho_t|_{t=0}=\vec v_0$. Since $\rho$ is irreducible, see \cite[Lemma 2.17]{BCLS}, we may choose 
$\{x_1,\ldots, x_{d+1}\}$ and $\{y_1,\ldots, y_{d+1}\}$ so that their images
$F=\{\xi_\rho^1(x_1),\dots, \xi_\rho^{1}(x_{d+1})\}$ and $F^*=\{\xi_\rho^{d-1}(y_1),\dots, \xi_\rho^{d-1}(y_{d+1})\}$ are projective frames.
Let $F_t=\{\xi_{\rho_t}^1(x_1),\dots, \xi_{\rho_t}^{1}(x_{d+1})\}$ and $F_t^*=\{\xi_{\rho_t}^{d-1}(y_1),\dots, \xi_{\rho_t}^{d-1}(y_{d+1})\}$. We may assume,
by restricting the path, that $F_t$ and $F_t^*$ are projective frames, and, by conjugating, that $F_t=F_0$ for all $t$.

Then  $\mu_{F_t,F_t^*}(\rho_t(\alpha))=[{\bf b}_{\rho_t}(y_i,z,\alpha(x_j),w)]$, for all $\alpha\in\Gamma$ and all $t$. So, since our cross-ratios
have log type zero, we see that 
$$\frac{d}{dt}\Big|_{t=0} \mu_{F_t,F_t^*}(\rho_t(\alpha))=0$$
for all $\alpha\in\Gamma$.

By construction, if $B\in\mathsf{PSL}(d,\mathbb R)$ and $F$ and $F^*$ are any projective frames, then $\mu_{F,B^*F^*}(A)=\mu_{F,F^*}(B^{-1}A)$.
So, if we choose $C_t\in\mathsf{PSL}(d,\mathbb R)$ so that $(C_t^{-1})^*F_t^*=F_0^*$,
then
$$0=\frac{d}{dt}\Big|_{t=0} \mu_{F_t,F_t^*}(\rho_t(\alpha))=\frac{d}{dt}\Big|_{t=0} \mu_{F_0,F_0^*}(C_t\rho_t(\alpha))=D\mu_{F_0,F_0^*} \left(\frac{d}{dt}\Big|_{t=0} C_t\rho_t(\alpha)\right)$$
for all $\alpha\in\Gamma$. Since $\mu_{F_0,F_0^*}$ is an immersion, this implies that
$$\frac{d}{dt}\Big|_{t=0} C_t\rho_t(\alpha)=0$$
for all $\alpha\in\Gamma$, so
$$C_0\circ \frac{d}{dt}\Big|_{t=0}\rho_t(\alpha)+\left(\frac{d}{dt}\Big|_{t=0} C_t\right)\circ \rho_t(\alpha)=0$$
for all $\alpha\in\Gamma$. By considering the case where $\alpha=id$, we see that  $\dot C_0=\left(\frac{d}{dt}\Big|_{t=0} C_t\right)=0$.
Since $C_0=I$, we see that $\frac{d}{dt}\Big|_{t=0}\rho_t(\alpha)=0$ for all $\alpha\in\Gamma$, which implies that $\vec v_0=0$.
\end{proof}

\section{The Hilbert length pressure metric}

Theorem \ref{BCLS generalized} has the following  immediate corollary:
\begin{corollary}
\label{stratum nondegen}
If $\mathsf{S}$ is a simple subgroup of $\mathsf{PSL}(d,\mathbb R)$ and 
$W$ is a submanifold of $\mathcal H_d(\Gamma)$ consisting of representations whose Zariski closure is $\mathsf{S}$,
then $\mathbb P^{\omega_H}$ is non-degenerate on $\TT W$.
\end{corollary}

\begin{proof}
Consider the Adjoint representation $Ad:\mathsf{S}\to\mathsf{SL}(V)$ where $V$ is the Lie algebra of $\mathsf S$.
Then $\mathsf{H}=Ad(\mathsf{S})$ is an irreducible reductive
subgroup of $\mathsf{SL}(V)$. 
Moreover, if $\rho\in W$, then
$Ad\circ \rho$ is irreducible and $\mathsf{H}$-generic. Theorem \ref{BCLS generalized} implies that the pressure form $\mathbb P^{\omega_1}$
is non-degenerate on $\TT Ad(W)$.

Note that $Ad$ is an immersion as the adjoint representation $ad\colon V\to\frak{sl}(V)$ is injective.
Recall that $\omega_1(Ad\circ\rho(\gamma))=\omega_H(\rho(\gamma))$. Therefore, $\mathbb P^{\omega_H}|_{\TT W}$ is the pull-back
of $\mathbb P^{\omega_1}|_{\TT Ad(W)}$ and hence non-degenerate.\end{proof}

The proof of \cite[Lemma 13.1]{BCLS} immediately generalizes to give the following lemma. (See Appendix \ref{app:path-metric} for a proof.)

\begin{lemma}
\label{pro:path-metric}
Let $W_0$ be a smooth manifold and let $W_n\subset W_{n-1}\subset \cdots W_1\subset W_0$ be a nested collection of  submanifolds of $W_0$
so that $W_i$ has positive codimension in $W_{i-1}$ for all $i$. Suppose that 
$g$ is a smooth non-negative symmetric 2-tensor $g$ such that 
\begin{itemize}
\item $g$ is positive definite on $\ms T_xW_{i-1}$ if $x\in W_{i-1}\setminus W_i$,
\item  the restriction of $g$ to $\ms T_xW_n$ is positive definite if $x\in W_n$.
\end{itemize}
Then the path pseudo metric defined by $g$ is a metric.
\end{lemma}

One thus obtains a pressure path metric on $\mathcal H_d(\Gamma)$ associated to the Hilbert  length functional, by applying
Sambarino's analysis of the possible Zariski closures of Hitchin representations (see Theorem \ref{thm: Sambarino}).

\medskip\noindent
{\bf Theorem \ref{hilbert mostly nondegenerate}.} {\em
Suppose that $\Gamma\subset\mathsf{PSL}(2,\mathbb R)$ is torsion-free and geometrically finite. If $\vec v\in \TT\mathcal H_d(\Gamma)$ is non-zero, then 
$\mathbb P^{\omega_H}(\vec v,\vec v)=0$ if and only if $\vec v$ is anti-self-dual. 
Moreover,
$\mathbb P^{\omega_H}$ gives rise to a mapping class group invariant path metric on $\mathcal H_d(\Gamma)$ which is an analytic Riemannian metric off of the self-dual locus.
}

\begin{proof} Let $W_0=\mathcal H_d(\Gamma)$. If $n\ge 4$ is even, let $W_1\subset\mathcal H_d(\Gamma)$ be the
submanifold of representations whose Zariski closure is conjugate into $\mathsf{PSp}(n,\mathbb R)$ and let $W_2$ be the Fuchsian locus (i.e. representations
contained in an irreducible image of $\mathsf{PSL}(2,\mathbb R)$.) If $n\ge 5$ is odd, and $n\ne 7$, let $W_1\subset\mathcal H_d(\Gamma)$ be the
submanifold of representations whose Zariski closure is conjugate into $\mathsf{PSO}(n,n-1)$. If $n=7$, let $W_1\subset\mathcal H_d(\Gamma)$ be the
submanifold of representations whose Zariski closure is conjugate into $\mathsf{PSO}(4,3)$, let $W_2\subset\mathcal H_d(\Gamma)$ be the
submanifold of $W_1$ consisting representations whose Zariski closure is conjugate into $\mathsf{G}_2$ and let $W_3$ be the Fuchsian locus.
If $n=3$, let $W_1$ be the Fuchsian locus. (Theorem \ref{FGcell} implies that $W_0$ is a manifold and that $W_i$ is
always a submanifold of $W_{i-1}$).

In all cases, Corollary \ref{stratum nondegen} implies that $\mathbb P^{\omega_H}$ is non-degenerate on $\ms T_xW_{i-1}$ if $x\in W_{i-1}\setminus W_i$.
Therefore, Lemma \ref{pro:path-metric}  implies that $\mathbb P^{\omega_H}$ gives rise to a path metric on $\mathcal H_d(\Gamma)$ which is Riemannian
off of $W_1$.
\end{proof}

\section{The (first) simple root pressure metric}\label{simple pressure}

\subsection{Trace functions} Recall that an element in $\mathcal H_d(\Gamma)$ is a conjugacy class of Hitchin representations of $\Gamma$ into $\mathsf{PSL}(d,\mathbb R)$. It will be convenient to identify $\mathcal H_d(\Gamma)$ with a subset $\widehat{\mathcal H}_d(\Gamma)$ of the character variety
$$X_d(\Gamma)=\mathrm{Hom}(\Gamma,\mathsf{SL}(d,\mathbb C))//\mathsf{SL}(d,\mathbb C).$$
The character variety $X_d(\Gamma)$ is the biggest Hausdorff quotient of $\mathrm{Hom}(\Gamma,\mathsf{SL}(d,\mathbb C))$ by the $\mathsf{SL}(d,\mathbb C)$-action by conjugation which coincides with the GIT quotient of $\mathrm{Hom}(\Gamma,\mathsf{SL}(d,\mathbb C))$ by this same action (See \cite{MFK git}).

If $\Gamma$ is cocompact, then Hitchin \cite{hitchin} showed that there is a  component, $\widehat{\mathcal H}_d(\Gamma)$, of $X_d(\Gamma)$ and
an analytic diffeomorphism $F:\mathcal H_d(\Gamma)\to\widehat{\mathcal H}_d(\Gamma)$, so that $F([\rho])$ is the conjugacy class of a lift to $\mathsf{SL}(d,\mathbb R)$ of $\rho$. 

Let $\widetilde{\mathcal H}_d(\Gamma)$ denote the set of all Hitchin representations of $\Gamma$ in $\mathsf{PSL}(d,\mathbb R)$.
If $\Gamma$ is not cocompact, then,
since $\Gamma$ is a free group and $\widetilde{\mathcal H}_d(\Gamma)$ is an analytic manifold, it is easy to define an analytic map 
$$F:\widetilde{\mathcal H}_d(\Gamma)\to \mathrm{Hom}(\Gamma,\mathsf{SL}(d,\mathbb C))$$
so that $F(\rho)$ is a lift of $\rho$. Since Hitchin representations are strongly irreducible, see Theorem \ref{cusped Hitchin properties},
Schur's Lemma implies that $F(\rho)$ is conjugate to $F(\eta)$ in $\mathsf{SL}(d,\mathbb C)$ if and only $\rho$ and $\eta$ are conjugate in $\mathsf{PSL}(d,\mathbb R)$.
Then, again since Hitchin representations are strongly irreducible, 
it follows that $F$ descends to an analytic embedding  $\hat F:\mathcal H_d(\Gamma)\to X_d(\Gamma)$ whose image lies in the smooth part of $X_d(\Gamma)$. 
We then let $\widehat{\mathcal H}_d(\Gamma)=\hat F(\mathcal H_d(\Gamma))$. Notice that if $d$ is odd, then $\widehat{\mathcal H}_d(\Gamma)=\mathcal H_d(\Gamma)$.
\medskip

If $\gamma\in\Gamma$, there is a complex analytic  trace function $\tr_\gamma:X_d(\Gamma)\to\mathbb C$ so that $\tr_\gamma([\rho])$ is the trace of $\rho(\gamma)$.
It is well-known that derivatives of trace functions generate the co-tangent space at any smooth point, see for example  Lubotzky-Magid \cite{lubotzky-magid}.
The following consequence will be used to verify the non-degeneracy of the first simple root pressure form.

\begin{lemma}
\label{traces generate}
If $[\rho]\in \widehat{\mathcal H}_d(\Gamma)$, then $\{{\rm D}_{\vec v}\tr_\g\mid\g\in\Gamma\}$ spans the cotangent space 
$\TT^*_{[\rho]}\widehat{\mathcal H}_d(\Gamma)$.
\end{lemma}

Even though $\tr_\gamma$ is not well-defined on $\mathcal H_d(\Gamma)$, we will abuse notation by saying that ${\rm D}_{\vec v}\mathrm{Tr}_\beta=0$
for some $\vec v\in \TT\mathcal H_d(S)$ if ${\rm D}_{D F(\vec v)}\mathrm{Tr}_\beta=0$.
\subsection{Nondegeneracy}

Bridgeman, Canary, Labourie and Sambarino \cite{BCLS2} prove that if $\Gamma$ is a closed surface group, then $\mathbb P^{\alpha_1}$ is non-degenerate
on $\mathcal H_d(\Gamma)$. A key tool in their work is the fact, due to Potrie-Sambarino \cite{potrie-sambarino}, that the topological entropy $h^{\alpha_1}(\rho)=1$ if
$\rho\in \mathcal H_d(\Gamma)$. Canary, Zhang and Zimmer \cite{CZZ2} generalized Potrie and Sambarino's result to the setting of torsion-free lattices which
are not cocompact.

\begin{theorem}[Potrie-Sambarino \cite{potrie-sambarino} and Canary-Zhang-Zimmer \cite{CZZ2}]
\label{root entropy 1}
If \hbox{$\Gamma\subset\mathsf{PSL}(2,\mathbb R)$}  is a torsion-free lattice, and
$\rho\in\mathcal H_d(\Gamma)$, then $h^{\alpha_1}(\rho)=1$.
\end{theorem}

With this result in hand, we are ready to establish the non-degeneracy of the first simple root pressure metric.

\medskip\noindent
{\bf Theorem \ref{simple nondegenerate}.} {\em
If $\Gamma\subset\mathsf{PSL}(2,\mathbb R)$ is a torsion-free lattice, then the pressure form ${\mathbb P}^{\alpha_1}$ is non-degenerate, so it
gives rise to a mapping class group invariant, analytic Riemannian metric on $\mathcal H_d(\Gamma)$.}

\medskip

Proposition \ref{pressure form basics} and Lemma \ref{traces generate}
together imply that Theorem \ref{simple nondegenerate} follows from the following proposition.

\begin{proposition}
\label{type zero translation}
If $\vec v\in T_{[\eta]}\mathcal H_d(\Gamma)$  and $\mathbb P^{\alpha_1}(\vec v,\vec v)=0$, then
${\rm D}_{\vec v}\tr_\beta=0$ for all $\beta\in\Gamma$.
\end{proposition}

Here, we will only sketch the proof, since the proof proceeds exactly as in the proof of \cite[Prop. 7.4]{BCLS2}. 

\begin{proof}
We again  abuse notation by identifying $[\rho]$ with  $F([\rho])\in\widehat{\mathcal H}_d(\Gamma)$.
Since $h^{\alpha_1}(\rho)=1$ for all $\rho\in\mathcal H_d(\Gamma)$, Proposition \ref{degeneracy condition} implies that ${\rm D}_{\vec v}\ell^{\alpha_1}_\beta=0$ for
all $\beta\in\Gamma$.

If $\alpha\in\Gamma$ is parabolic, then $\tr_\alpha$  is constant on $\widehat{\mathcal H}_d(\Gamma)$, so ${\rm D}_{\vec v}\tr_\alpha=0$.
 
If $\beta$ is hyperbolic,  we may choose $\alpha\in\Gamma$, so that $\alpha$ is hyperbolic and $\alpha$ and $\beta$ have non-intersecting axes.
We may pass to powers  $\alpha^n$ and $\beta^n$ which generate a Schottky subgroup of $\Gamma$. 
We are then exactly in the setting of the proof of \cite[Prop. 7.4]{BCLS2} which shows that 
${\rm D}_{\vec v}\lambda_i(\rho(\beta^n))= {\rm D}_{\vec v}\lambda_i(\rho(\beta))^n=0$ for all $i$. Therefore, 
${\rm D}_{\vec v}\lambda_i(\rho(\beta))= 0$ for all $i$, so ${\rm D}_{\vec v}\tr_\beta= 0$.
\end{proof}

\section{Appendix}\label{app:path-metric}

We prove:

\medskip\noindent
{\bf Lemma \ref{pro:path-metric}.}
{\em
Let $W_0$ be a smooth manifold and let $W_n\subset W_{n-1}\subset \cdots W_1\subset W_0$ be a nested collections of submanifolds of $W_0$
so that $W_i$ has non-zero codimension in $W_{i-1}$ for all $i$. Set $W_{n+1}=\emptyset$. Suppose that 
$g$ is a smooth non-negative symmetric 2-tensor on $W_0$ such that for every $i=0,\dots, n$, the restriction of $g$ to $\ms T_xW_i$ is positive definite if $x\in W_i\setminus W_{i+1}$. Then, the path pseudo-metric defined by $g$ is a metric.}

\begin{proof}
We proceed iteratively to establish the following claim:

\medskip\noindent
{\bf Claim.} {\em If $x\in W_i\setminus W_{i+1}$, then $x$ has a neighborhood $U$ whose closure $\overline U$ lies in $W_0\setminus W_{i+1}$, so that if $u\in \overline U\setminus\{x\}$, then $d(x,u)>0$.}
\medskip

Once we have proved this claim for all $i$, we will have completed the proof of the lemma.

If $x\in W_0\setminus W_1$, then if $U$ is any neighborhood of $x$ whose closure $\overline U$ is disjoint from $W_1$, then $g$ is Riemannian on $\overline U$.
Therefore, our claim is true for $i=0$. 

Next, we suppose that the claim is true for all $i=j<k$, and prove the claim for $i=k$. This establishes the claim for all $i$.

Let $n_i=\dim W_i$. If $x\in W_k\setminus W_{k+1}$,
we may identify some neighborhood  $U$ of $x$ with the Euclidean unit ball  in $\mathbb R^{n_0}$ (centered at $\vec 0$) so that $x$ is identified with $\vec 0$. 
We may assume that the closure $\overline U$ of $U$ is compact and disjoint from $W_{k+1}$ and
that if $j\le k$, then $W_j\cap \overline U$ is identified with the intersection of the closure $D(\vec 0,1)$ of $B(\vec 0,1)$ with $\mathbb R^{n_j}\times \{\vec 0\}^{n_0-n_j}$. 
We will work in coordinates for the rest of this proof.
We identify $TD(\vec 0,1)$ with $D(\vec 0,1)\times \mathbb R^{n_0}$.

Since the restriction of $g$ to $T(W_k\setminus W_{k+1})$ is Riemannian, there exists $r,s>0$ so that if $\vec v$ is a (Euclidean) unit vector in
$(W_k\cap D(\vec 0,1))\times \mathbb R^{n_k}\times \{\vec 0\}^{n_0- n_k}$, then $s^2\ge g(\vec v,\vec v)\ge 4r^2$. 

Since $g$ is continuous, it follows that, after possibly restricting to a smaller neighborhood of $x$, we can assume
that if $\vec v$ is a unit vector in  $D(\vec 0,1)\times \mathbb R^{n_k}\times \{\vec 0\}^{n_0- n_k}$, then $4s^2\ge g(\vec v,\vec v)\ge r^2$. It follows that the (Euclidean) projection map
from $\pi_k\colon D(\vec 0,1)\to W_k$ is $K$-Lipschitz where $K=\frac{2s}{r}$. Therefore, since the restriction of $g$ to $T(W_k\setminus W_{k+1})$ is Riemannian, it follows that if
$u\in U$ and $\pi_k(u)\ne \vec 0$, then $d(u,x)>0$. On the other hand, if $\pi_k(u)=\vec 0$ and $u\ne x$, then $u\in W_0\setminus W_k$, so, by our iterative assumption,
there exists a neighborhood $V$ of $u$ whose closure lies in $W_0\setminus W_k$, so that if $v\in \overline V\setminus\{u\}$, then $d(v,u)>0$. It follows that there exists $c>0$ so that if $v\in\partial V$, then $d(u,v)\ge c$.
Since $x\notin V$, this implies that $d(x,u)\ge d(\partial V, u)\ge c>0$. This completes the proof of the claim and hence the lemma.
\end{proof}

\end{document}